\def\k{{\kappa}}

\def\p{{\pi}}

\def\cB{{\cal B}}
\def\cC{{\cal C}}

\def\cG{{\cal G}}

\def\zN{{\mathbb N}}

\def\pf{{\hfill$\Box$}}
\def\pfc{{\hfill$\diamondsuit$}}

\def\rar{{\rightarrow}}

\def\diam{{\rm diam}}
\def\proof{{\noindent {\it Proof}.\ \ }}

\def\l({{\left(}}
\def\r){{\right)}}

\def\({{\Biggl(}}
\def\){{\Biggr)}}
\def\[{{\Biggl[}}
\def\]{{\Biggr]}}

\documentclass[12pt]{article}

\newtheorem{thm}{Theorem}
\newtheorem{clm}[thm]{Claim}
\newtheorem{cnj}[thm]{Conjecture}
\newtheorem{cor}[thm]{Corollary}

\newtheorem{lem}[thm]{Lemma}
\newtheorem{prb}[thm]{Problem}
\newtheorem{prp}[thm]{Proposition}
\newtheorem{qst}[thm]{Question}
\newtheorem{res}[thm]{Result}

\usepackage{graphicx,amssymb,amsmath,verbatim}

\begin{document}

%
%
\title{On Pebbling Graphs by their Blocks}

\author{
Dawn Curtis\thanks{dawn.curtis@asu.edu},
Taylor Hines\thanks{taylor.hines@asu.edu},
Glenn Hurlbert\thanks{hurlbert@asu.edu},
Tatiana Moyer\thanks{tatiana.moyer@asu.edu}\\
    Department of Mathematics and Statistics\\
    Arizona State University,
    Tempe, AZ 85287-1804
}
\maketitle

\newpage


%
%
\begin{abstract}

Graph pebbling is a game played on a connected graph $G$.
A player purchases pebbles at a dollar a piece, and hands them to an
adversary who distributes them among the vertices of $G$ (called a
configuration) and chooses a target vertex $r$.
The player may make a pebbling move by taking two pebbles off of one vertex
and moving one pebble to a neighboring vertex.
The player wins the game if he can move $k$ pebbles to $r$.
The value of the game $(G,k)$, called the
$k$-pebbling number of $G$ and denoted $\p_k(G)$, is the minimum cost
to the player to guarantee a win.
That is, it is the smallest positive integer $m$ of pebbles so that,
from every configuration of size $m$, one can move $k$ pebbles to any target.
In this paper, we use the block structure of graphs to investigate pebbling
numbers, and we present the exact pebbling number of the graphs whose blocks
are complete.
We also provide an upper bound for the $k$-pebbling number of diameter-two
graphs, which can be the basis for further investigation into the pebbling numbers
of graphs with blocks that have diameter at most two.

\end{abstract}

\newpage

%
%
\section{Introduction}\label{Intro}

Graph pebbling is a game played on a connected graph $G=(V,E)$.%
\footnote{We assume the notation and terminology of \cite{W} throughout.}
A player purchases pebbles at a dollar a piece, and hands them to an
adversary who distributes them among the vertices of $G$ (called a
{\it configuration}) and chooses a target, or {\it root} vertex $r$.
The player may make a {\it pebbling move} by taking two pebbles off of one vertex
and moving one pebble to a neighboring vertex.
The player wins the game if he can move $k$ pebbles to $r$, in which case
we say that $r$ is $k$-{\it pebbled}.
Another common terminology calls the configuration $k$-{\it fold} $r$-{\it solvable}.
The {\it value} of the game $(G,k)$, called the
$k$-{\it pebbling number} of $G$ and denoted $\p_k(G)$, is the minimum cost
to the player to guarantee a win.
That is, it is the smallest positive integer $m$ of pebbles so that,
from every configuration of size $m$, one can move $k$ pebbles to any root.
If $k$ is not specified, it is assumed to be one.

For example, by the pigeonhole principle we have $\p(K_n)=n$, where $K_n$
is the complete graph on $n$ vertices.
From there, induction shows that $\p_k(K_n)=n+2(k-1)$.
Induction also proves that $\p_k(P_n)=k2^{n-1}$, where $P_n$ is the path
on $n$ vertices.
These two graphs illustrate the tightness of the two main lower bounds
$\p(G)\ge\max\{n(G),2^{\diam(G)}\}$, where $\diam(G)$ is the {\it diameter}
of $G$, the number of edges in a maximum induced path.
Another fundamental result uses the path fact and induction to calculate
the $k$-pebbling number of trees (see \cite{Chung}).
The survey \cite{H} contains a wealth of information regarding pebbling
results and variations.

Complete graphs and paths are examples of {\it greedy} graphs.
That is, the most efficient pebbling moves are directed towards the root.
More formally, a pebbling move from $u$ to $v$ is greedy if
$dist(v,r) < dist(u,r)$, where $dist(x,y)$ denotes the distance between $x$
and $y$.
A greedy solution uses only greedy moves.
A graph $G$ is greedy if every configuration of size $\p(G)$ can be greedily
solved.
If a graph is greedy, then we can assume every pebbling move is directed
towards the root.
The greedy property of trees follows from the No-Cycle Lemma of \cite{Moe} (see also \cite{CCF,MC}), which states that the digraph whose arcs represent the
pebbling moves of a minimal solution contains no directed cycles.
A {\it cut vertex} of a graph is a vertex that, if removed, disconnects the graph.
The {\it connectivity} $\k$ of a graph is the minimum number of
vertices whose deletion disconnects the graph or reduces it to only one vertex.
Two important results relate diameter and connectivity to pebbling numbers.
Pachter, Snevily, and Voxman proved the first.

\begin{res}\label{psv}\cite{PSV}.
If $G$ is a connected graph on $n$ vertices with $\diam(G)\le 2$
then $\p(G)\le n+1$.
\end{res}

Clarke, Hochberg, and Hurlbert \cite{CHH} characterized which diameter
two graphs have pebbling number $n$ and which have pebbling number $n+1$.
We will use the graphs that describe that characterization in Section
\ref{Cliques}.
Motivated by the characterization, Czygrinow, Hurlbert, Kierstead, and Trotter
proved the second.

\begin{res}\label{chkt}\cite{CHKT}.
If $G$ is a connected graph on $n$ vertices with $\diam(G)\le d$
and $\k(G)\ge 2^{2d+3}$ then $\p(G)=n$.
\end{res}

\begin{figure}
\centerline{\includegraphics[height=3.0 in]{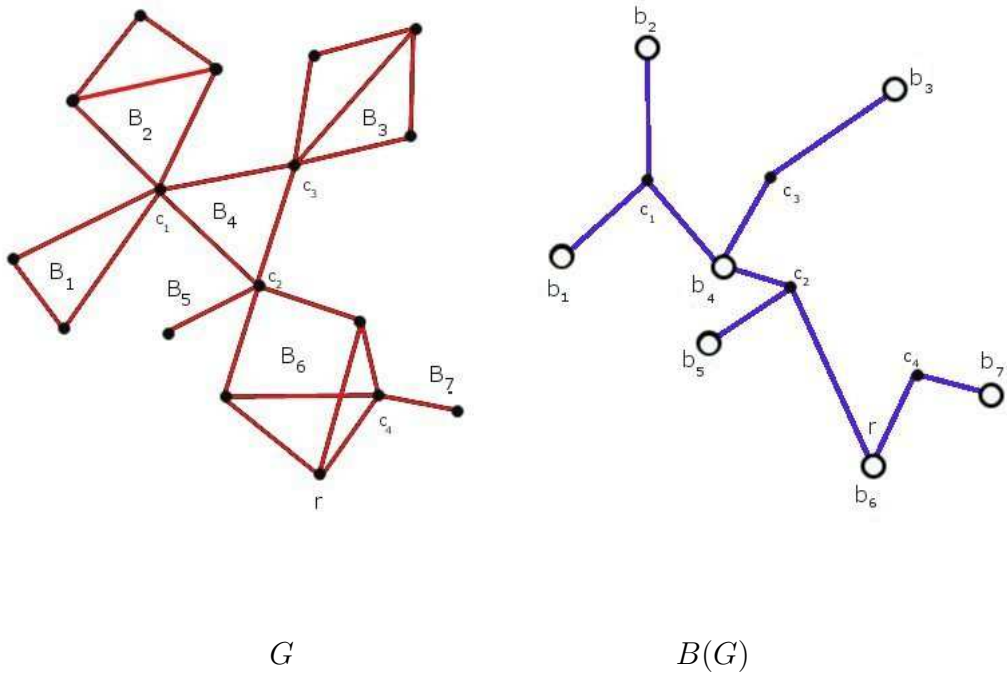}}
$$G\hspace{2.0 in}B(G)$$
\caption{A graph and its block-cutpoint graph}
\label{BofG}
\end{figure}

This result states that high connectivity compensates for large
diameter in keeping the pebbling number to a minimum.
In this paper we exploit graph structures further to investigate
pebbling numbers.
A {\it block} of a graph $G$ is a maximal subgraph of $G$ with no cut
vertex.
Let $\cB$ be the set of all blocks of $G$ and $\cC$ be the set of
all cut vertices of $G$.
Then the {\it block-cutpoint} graph of $G$, denoted $B(G)$, has
vertices $\cB\cup\cC$, with edges $(B,C)$ whenever $C\in V(B)$.
Note that $B(G)$ is always a tree (see \cite{W}).
Figure \ref{BofG} shows an example.

Here we instigate a line of research into using the $k$-pebbling numbers
of $B(G)$ and of the blocks of $G$ to give upper bounds on $\p_k(G)$.
To begin, we generalize Chung's tree result to weighted trees in Section
\ref{Trees}.
We then present the exact $k$-pebbling number of $G$ when every block of
$G$ is complete in Section \ref{Cliques}.
Also in Section \ref{Cliques}, we prove the following theorem, and show that
there is a diameter-2 graph $G$ on $n\ge 6$ vertices with $\p_k(G)=n+4k-3$
for all $k$ (Theorem \ref{bigconfig}).
Thus Theorem \ref{kpd2} is not known to be tight.

\begin{thm}\label{kpd2}
If $G$ is a graph on $n$ vertices with $\diam(G)\le 2$
then $\p_k(G)\le n+7k-6$.
\end{thm}

Section \ref{Remarks} provides some further conjectures, questions, and
possibilities for future research.

%
%
\section{Trees and General Pebbling}\label{Trees}

\begin{figure}
\centerline{\includegraphics[height=3.0 in]{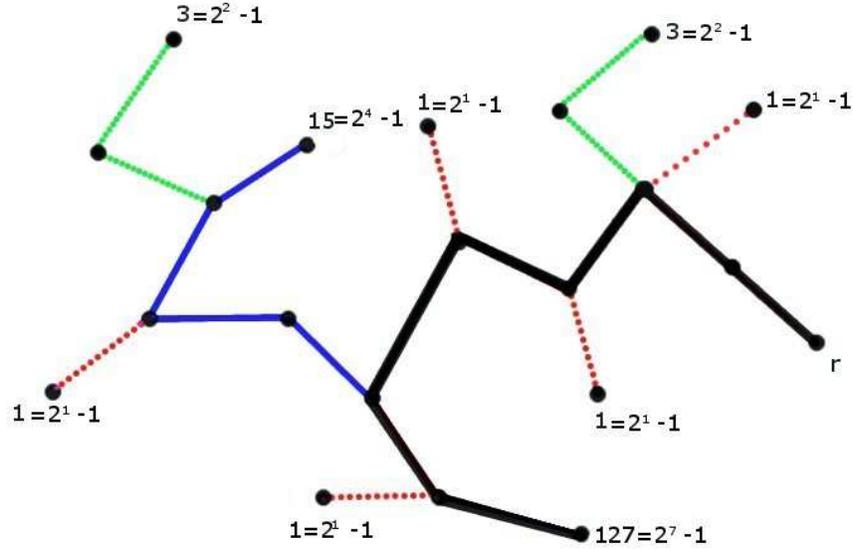}}
\caption{A non-maximal $r$-path partition of a tree, with its corresponding
unsolvable configuration}
\label{NotMaxTree}
\end{figure}

\begin{figure}
\centerline{\includegraphics[height=3.0 in]{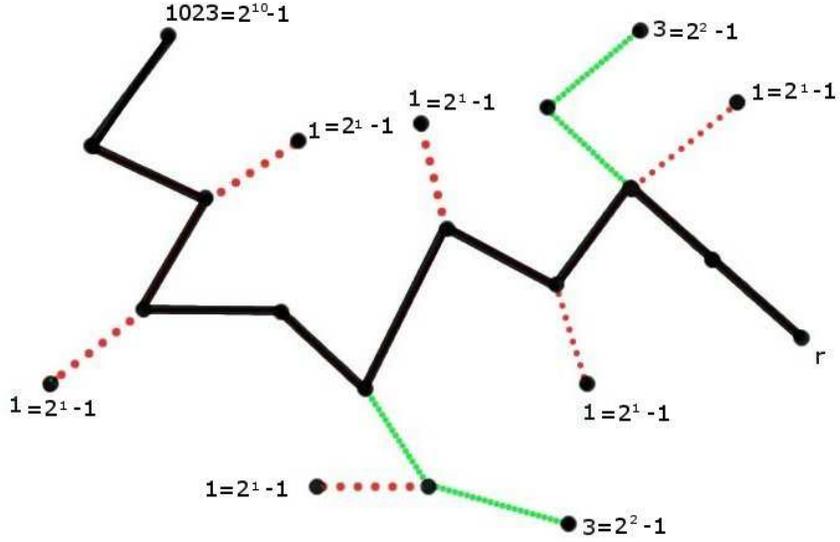}}
\caption{An $r$-maximal path partition of a tree, with its corresponding
unsolvable configuration}
\label{Max}
\end{figure}

A {\it tree} is a connected, acyclic graph, and a {\it forest} is a union of
pairwise vertex-disjoint trees.
A {\it leaf} of a tree is a vertex of degree one.
An $r$-{\it path partition} of a particular tree $T$ is a partition of the edges 
of $T$ into paths, constructed by carrying out the following algorithm.
Construct the sequence of pairs ${(T_i, F_i)}$, where each $T_i$ is a tree and each
$F_i$ is a forest, with $E(T_i) \cup E(F_i) = E(T)$, and $E(T_i) \cap E(F_i) = \emptyset$.
Begin with $T_0 = r$, $F_0 = T$ and end with $T_t = T$, and $F_t = \emptyset$.
At each stage, for some path $P_i$ we have $P_i = T_i - T_{i-1} = F_{i-1} - F_i$, with
the property that for each i, the intersection $V(P_i)\cap V(T_{i-1})$ 
is a leaf of $P_i$.
The path partition is {\it r-maximal} if each $P_i$ is the longest such path in
$F_{i-1}$.
An $r$-maximal path partition is {\it maximal} if $r$ is one of the leaves of the longest
path in $T$.
An r-path partition of a tree is depicted in Figure \ref{NotMaxTree}, and a maximal
path partition of a tree is depicted in Figure \ref{Max}.

Define $x_i$ to be the leaf of $P_i$ in $T_{i-1}$ and $y_i$ to be the leaf of $P_i$ not
in $T_{i-1}$, and let $a_i = |E(P_i)|$.

\begin{lem}\label{Unsolvable}.
The configuration $C$ on $T$ defined by each $C(y_i)=2^{a_i}-1$ and $C(v)=0$
for all other $v$ is $r$-unsolvable.
\end{lem}

\proof
We use induction.
Let $C_i$ be the restriction of $C$ to $T_i$.
The case in which $i=0$ is trivial since the root has no pebbles.
Now, assume that $C_k$ is $r$-unsolvable on $T_k$.
We know that the configuration on $P_{k+1}$ is $x_{k+1}$-unsolvable because the
pebbling number of a path of length $l$ is $2^l$.
Thus, no pebbles can be moved to from $P_{k+1}$ to $T_k$ since
$V(T_k) \cap V(T_{k+1}) = x_{k+1}$.
Since we already know $T_k$ is unsolvable, $T_{k+1}$ is unsolvable also.
Thus, by induction, the configuration $C$ on $T$ is $r$-unsolvable.
\pf
\medskip

Chung's result generalizes this idea for k-pebbling.

\begin{res}\label{chung}\cite{Chung}.
If $T$ is a tree and $a_1, a_2,\ldots,a_t$ is the sequence of the path size (i.e.
the number of vertices in the path) in a maximum path partition of $T$, then
$\p_k(T)=k2^{a_1}+ \sum_{i=2}^{t} 2^{a_i} - t + 1$.
\end{res}

Chung's proof of this result uses induction performed on the vertices of $T$ by
fixing and then removing the root, thus dividing $T$ into subtrees in order to
use induction.
We give a different proof of the more general Theorem \ref{ourtree}, relying on
the fact that trees are greedy.

First we consider a more general form of pebbling.
For each edge $e$ of a graph $G$ we can assign a weight $w_e$.
The weight is intended to signify that it takes $w_e$ pebbles at one
end of $e$ to place 1 pebble at its other end.
Hence the pebbling considered to this point has $w_e=2$ for all $e$.
We define the weighted pebbling number $\p_k^w(G,r)$ to be the minimum number 
$m$ so that every configuration of size $m$ can $k$-pebble $r$ by using $w$-weighted 
pebbling moves on $G$.

\begin{figure}
\centerline{\includegraphics[height=3.0 in]{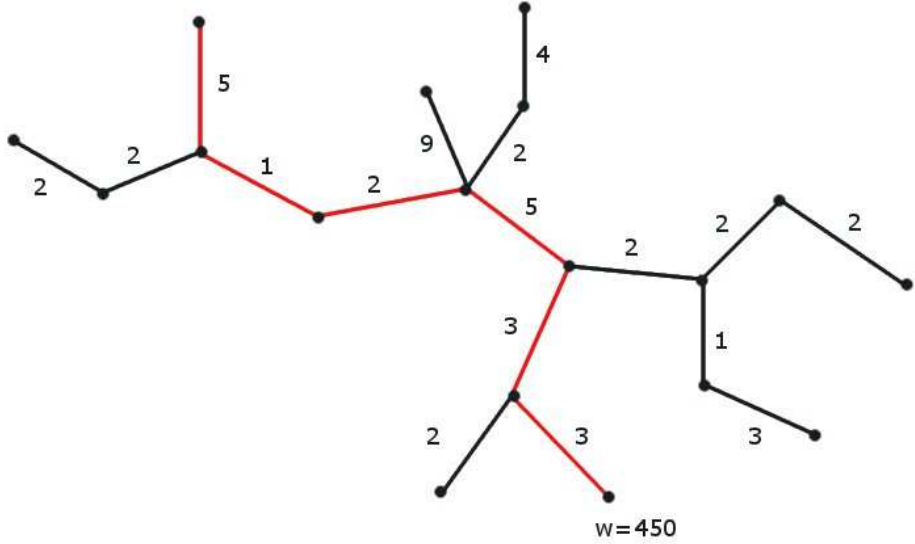}}
\caption{An edge-weighted tree}
\label{Weighted}
\end{figure}

Given a weight function $w:E(G)\rar\zN$, we extrapolate to a weight
function on the set of all paths of $G$, where $w(P)$ is the product
of edge weights over all edges of the path $P$.
Now when constructing maximal path partitions, we replace the condition
``longest path" by ``heaviest path" (greatest weight).
This is equivalent for constant weight 2 pebbling.
Nothing in the proof of Chung's theorem changes for weighted trees, but we
introduce a new proof of the pebbling number of a weighted tree.

Let $P_1,\ldots,P_t$ be an $r$-maximal path partition of $T$,
with $w(P_1) \ge \ldots \ge w(P_t)$.
Let $f_k^w(T,r)=kw(P_1)+ \sum_{i=2}^t w(P_i) - t + 1$.
For vertices $x$ and $y$ on a path $P$, denote by $P[x,y]$ the subpath of $P$
from $x$ to $y$.

\begin{thm}\label{ourtree}
Every weighted tree T satisfies $\p_k^w(T,r) \le f_k^w(T,r)$.
\end{thm}

\proof
The theorem is trivially true when $t=1$ since $T$ is a path.

For $t \ge 1$, define $T'=T-P_t$.
Then $f_k(T,r)=f_k(T',r)+w(P_t)-1.$
Let $P_j$ be a path containing the non-leaf endpoint $x_t$ of $P_t$,
and let vertex $y_j$ be the leaf of $T$ on $P_j$.
Define $W = w(P_j[x_t,y_j])$.
Thus we know from the maximal $r$-path construction that $W \ge w(P_t)$.

Let $C$ be an unsolvable configuration on $T$ with $|C|=f_k(T,r)$.
Without loss of generality, we can assume that all the pebbles are on the
leaves of a tree because the maximum sized unsolvable configuration sits
on the leaves only.
Let $s \ge 0$ be the number of pebbles $P_t$ contributes to the vertex $x_t$,
so we have $sw(P_t) \le |C(P_t)| < (s+1)w(P_t)$.

Now define the configuration $C'$ on $T'$ by $C'(y_j)=C(y_j)+sW$ and
$C'(v)=C(v)$ otherwise.
Then,
\begin{eqnarray*}
|C'|& = &|C|-[(s+1)w(P_t)-1] + sW\\
&\ge& f_k(T,r)-w(P_t)+1\\
&=& f_k(T',r)\ .
\end{eqnarray*}
Hence $C'$ is $k$-fold solvable on $T'$.
Now define $C^*$ on $T'$ by $C^*(x_t)=C(x_t)+s$ and $C^*(v)=C(v)$ otherwise.
In particular, because of greediness,
$C^*$ is $k$-fold $r$-solvable on $T'$ because moving at most $sw(P_t)$
pebbles from $y_j$ to $x_t$ converts $C'$ to a solvable subconfiguration of $C^*$.
Now, since $C(P_t) \ge sw(P_t)$, the base case says we can move $s$ pebbles
from $P_t$ to $x_t$, and in doing so we arrive again at $C^*$ on $T'$.
Hence $C$ is $k$-fold $r$-solvable.
\pf
\medskip

We will use Theorem \ref{ourtree} to upper bound the pebbling number of graphs
composed of blocks.
The technique utilizes the block-cutpoint graph.

For a graph $G$ and its block-cutpoint graph $B(G)$, let $b_i$ denote the
vertex of $B(G)$ that corresponds to the block $B_i$ in $G$.
For each block $B_i$, let $x_i$ denote the cut vertex of $G$ in $B_i$ that is
closest to the root (it is possible that some $x_i=x_j$).
Let $e_i$ denote the edge of $B(G)$ between $b_i$ and $x_i$, and define its
weight by $w(e_i)=\p(B_i, x_i)$.
Let all other edges have weight 1.
For a root $r$ of $G$, let $B$ denote the block containing it, represented by
the vertex $b$ in $B(G)$.
Let $B'(G)$ be the graph obtained from $B(G)$ by adjoining to $b$ by an edge of
weight 1 a new vertex $r'$.
Then we arrive at the following theorem.

\begin{thm}\label{BlockCut}
Every graph $G$ satisfies $\p_k(G,r) \le \p_k^w(B'(G),r')$
\end{thm}

\proof
For a set $U$ of vertices, denote by $C(U)$ the sum $\sum_{v \in U} C(v)$.
Let $x(B_i)$ denote all the cut vertices of $G$ in the block $B_i$.
Given a configuration $C$ on $G$, define $C'$ on $B'(G)$ by
\begin{itemize}
\item
$C'(x_i) = C(x_i)$ for all cut vertices $x_i$, and
\item
$C'(b_i) = C(B_i) - C(x(B_i))$ for all blocks $B_i$.
\end{itemize}
Given an $r'$-solution $S'$ of $C'$ on $B'(G)$, which exists because
$|C'|=|C|=\p_k^w(B'(G),r')$, define the $r$-solution $S$ of $C$
on $G$ by the following: replace every pebbling step along $e_i$ in $S'$ by some
$x_i$-solution of some $\p(B_i)$ of the pebbles in $B_i$.
Then $S$ is an $r$-solution.
\pf

%
%
\section{Larger Blocks}\label{Cliques}

In this section we consider the cases in which all blocks are cliques or all
have bounded diameters.
The following proposition is well known.

\begin{prp}\label{super}
If $H$ is a connected spanning subgraph of $G$ then $\p_k(G,r)\le \p_k(H,r)$
for every root $r$.
\end{prp}

Proposition \ref{super} holds because $r$-solutions in $H$ are $r$-solutions in
$G$.
In particular, this holds when $H$ is a breadth-first search spanning tree of
$G$ that is rooted at $r$ and thus preserves distances to $r$ in $G$.
This allows us to prove the following.

\begin{res}\label{CliqBlok}
Let $G$ be a connected graph in which every block is a clique.
Let $T$ be a breadth-first search spanning tree of $G$.
Then $\p_k(G)=\p_k(T)$.
\end{res}

\begin{figure}
\centerline{\includegraphics[height=3.0 in]{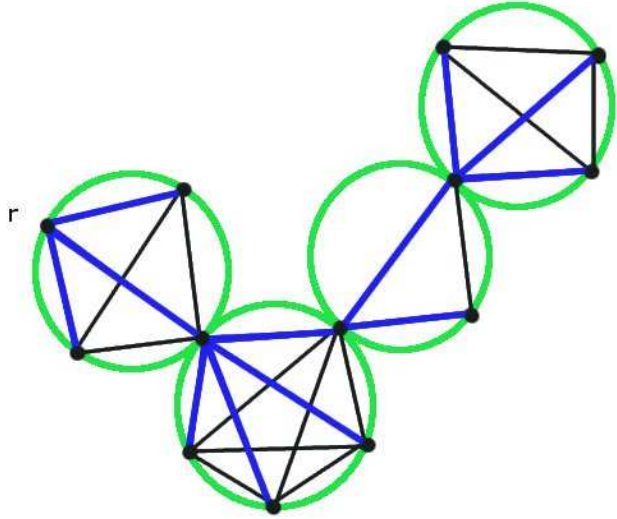}}
\caption{A clique block graph with its breadth-first search spanning tree}
\label{CBG}
\end{figure}

\proof
The fact that $\p_k(G) \le \p_k(T)$ follows from Proposition \ref{super}.
The fact that $\p_k(G) \ge \p_k(T)$ follows from showing that every
$r$-solvable configuration $C$ on $G$ is $r$-solvable on $T$.
Indeed, let $S$ be an $r$-solution in $G$, and for a block $B$ of $G$, denote
by $x=x(B)$ the cut vertex of $B$ that is closest to $r$.
If the sequence is greedy, then all its edges are in $T$.
If the sequence is not greedy, then $S$ contains an edge from some vertex $a$
to some vertex $b \neq x$.
Replace this edge by the edge from $a$ to $x$.
The resulting sequence is an $r$-solution on $T$.
Thus $\p_k(G) = \p_k(T)$.
\pf

\begin{cor}\label{Formula}
Let $G$ be a connected graph in which every block is a clique.
Let $T$ be a breadth-first search spanning tree of $G$.
Let ${a_1,\ldots,a_t}$ denote the path lengths in a maximal path partition of
$T$ rooted at $r$.
Then $\p_k(G,r)=n+2^{a_1}(k-1)+\sum_{i=1}^t(2^{a_i}-a_i-1)$.
\end{cor}

Note that the formula in Corollary \ref{Formula} is of the form $n+c_1k+c_2$,
which is also the form of the formula in Theorem \ref{kpd2}.
Also, the {\it fractional pebbling number}, defined as
$\hat{\p}(G)=\lim_{k\to\infty}\p_k(G)/k$ is seen to be $\hat{\p}(G)= 2^{diam(G)}$
for such $G$.
This is an instance of the Fractional Pebbling Conjecture of \cite{H}, recently
proven in \cite{HV}.

\begin{figure}
\centerline{\includegraphics[height=3.0 in]{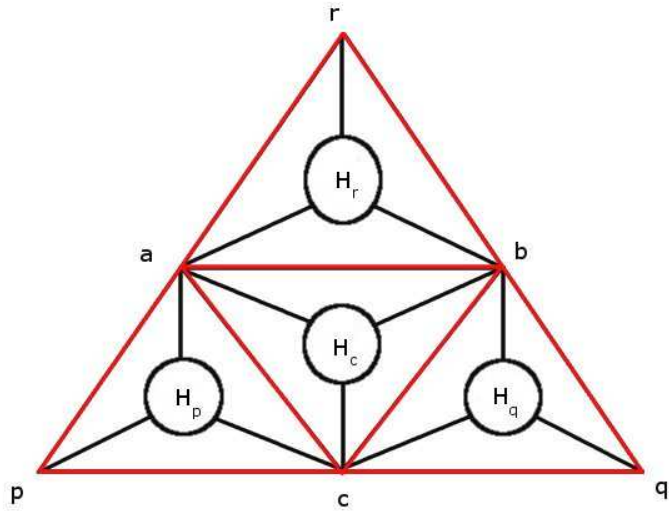}}
\caption{The extremal graph $\cal G$}
\label{Extremal}
\end{figure}

Now we provide the upper and lower bounds on diameter-two graphs.
To show a lower bound, we will display an unsolvable
configuration on an extremal graph $\cal G$.
This is the graph that Clarke, et al. \cite{CHH} used to
characterize the diameter two graphs with pebbling number $n+1$.
The vertices of $\cal G$ are $\{a,b,c,p,q,r\}\cup_{z\in\{p,q,r,c\}}V(H_z)$,
where $H_p,H_q,H_r$, and $H_c$ are any graphs with the following properties.
\begin{itemize}
\item
Every component of $H_p,H_q$, and $H_r$ has some vertex adjacent to $p,q$, and
$r$, respectively.
\item
Every vertex of $H_p,H_q$, and $H_r$ is adjacent to $a$ and $c$, $b$ and $c$,
$a$ and $b$, respectively.
\item
Every vertex of $H_c$ is adjacent to $a,b$, and $c$.
\end{itemize}
Furthermore, $(a,r,b,q,c,p)$ forms a 6-cycle, $(a,b,c)$ forms a triangle, as shown
in Figure \ref{Extremal}, and no other edges than previously mentioned are included.
Note that the diameter of $\cal G$ is 2.

\begin{thm}\label{bigconfig}
For all $n\ge 6$, there is a graph $G$ on $n$ vertices with
$\p_k(G) \ge n + 4k - 3$ for all $k$.
\end{thm}

\proof
As suggested above, we show that $\cal G$ is such a graph.
Distribute the following configuration of size $n + 4k - 4$ on the
$\cG$:
\begin{itemize}
\item
Place $4k - 1$ pebbles on $p$
\item
Place 3 pebbles on $q$
\item
Place 1 pebble on every vertex in $\cup_{z\in\{p,q,r,c\}}H_z$ and 0 elsewhere.
\end{itemize}
The configuration is $r$-unsolvable since every solution costs at least 4 pebbles
(because the pebbles are at distance 2 from $r$, and so after $k=1$ solutions
at most $n$ pebbles remain).
In fact, the remaining configuration is a subconfiguration of the one defined
above for $k=1$, which was shown to be $r$-unsolvable in \cite{PSV}.
Hence $\p_k(\cG) > n + 4k - 4$.
\pf
\medskip

\begin{figure}
\centerline{\includegraphics[width=5.0 in]{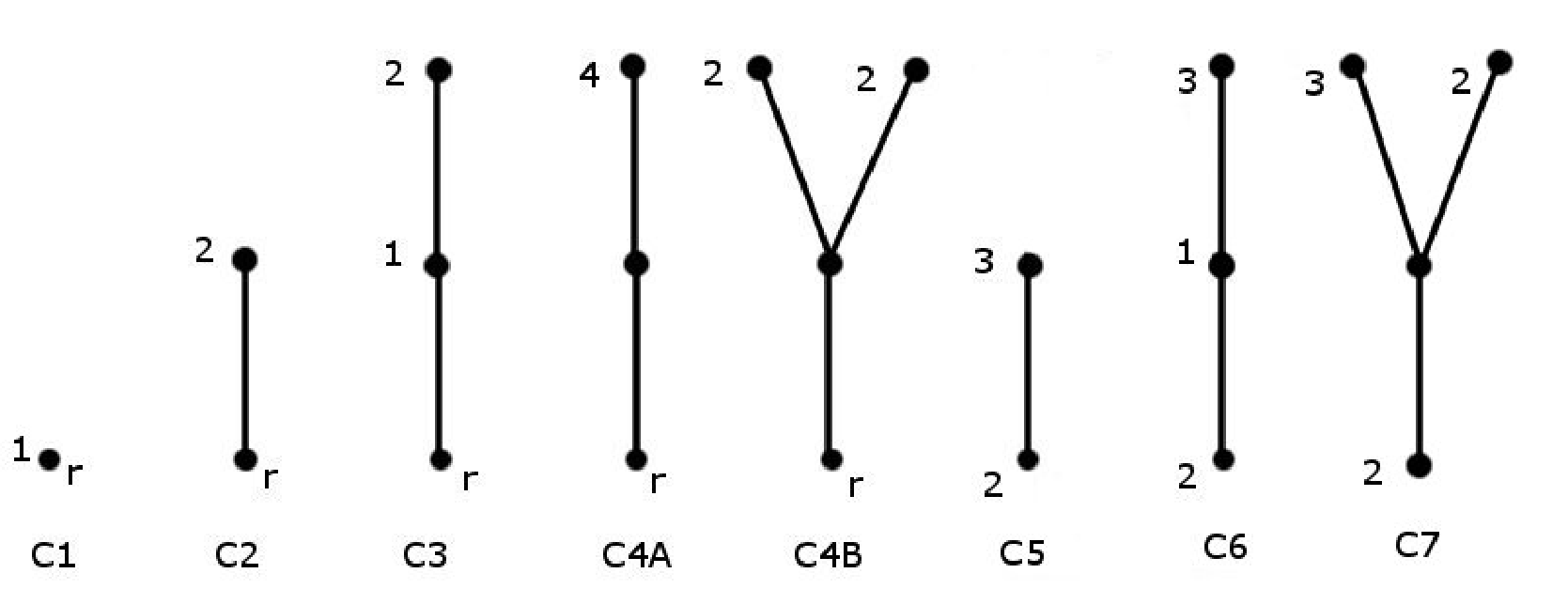}}
\caption{Cheap solutions of cost 7 or less}
\label{CheapSolns}
\end{figure}

To prove Theorem \ref{kpd2} we consider the eight
{\it cheap} configurations shown in Figure \ref{CheapSolns}.
We call them cheap because they lose a small number (at most 7)
of pebbles in the process of moving one pebble to the root.
In particular, their names indicate their {\it cost} (number of pebbles used).
For example, in {\sf C7}, {\sf C6}, and {\sf C5}, one moves an extra pebble onto
where 3 sits to create {\sf C4A}.
Then one can reach {\sf C2} from {\sf C4B}, {\sf C4A}, and {\sf C3}.
Of course, {\sf C2} results in {\sf C1}.
There are more cheap solutions than these, but we do not need them
in our argument.

We show by contradiction that a cheap solution must exist, and thus a pebble
can be moved to the root with the loss of at most 7 pebbles.
The remaining $k-1$ solutions will be found by induction.
\bigskip

\noindent {\it Proof of Theorem \ref{kpd2}.}
Assume that the configuration $C$ of pebbles on $G$ is of size
$n + 7k - 6$ and has no cheap solutions of cost $7$ or less.
We will derive a contradiction to show that a cheap solution exists.
Then after using a cheap solution we apply induction to get the
remaining $k-1$ solutions.
The theorem is already true for $k=1$ by Result \ref{psv}.
Define the following notation.
\begin{itemize}
\item
$N_i$ is the set of vertices with $i$ pebbles.
\item
$N_{i,r}$ is the set of common neighbors of $N_i$ and root $r$.
\item
$N_{i,j}$ is the set of common neighbors of pairs of vertices from $N_i$ and $N_j$.
\item
$n_i = |N_i|$, $n_{i,j}=|N_{i,j}|$, $n_{i,r}=|N_{i,r}|$, and $n_0'=|N_0'|$.
\item
$N_0'$ = $N_0 - N_{3,r} - N_{3,3} - N_{2,r}$.
\end{itemize}

\begin{clm}\label{claims}
If $C$ is a configuration on a diameter-2 graph $G$ with no cheap solutions, then
\begin{enumerate}
\item[{\sf S1.}]
$N_{i,r} \subseteq N_0$ for $i\in \{2,3\}$,
\item[{\sf S2.}]
$N_{3,3} \subseteq N_0$,
\item[{\sf S3.}]
$n_{i,r} \ge n_i$ for $i \in \{2,3\}$,
\item[{\sf S4.}]
$|C| = 3n_3 + 2n_2 + n_1$,
\item[{\sf S5.}]
$n = n_3 + n_2 + n_1 + (n_{3,r} + n_{3,3} + n_{2,r} + n_0')$, and
\item[{\sf S6.}]
$n_{3,3} \ge \binom{n_3}{2}$.
\end{enumerate}
\end{clm}

\noindent
{\it Proof of Claim \ref{claims}.}
We refer to Figure \ref{CheapSolns}.
Statement {\sf S1} follows from the nonexistence of {\sf C3} because a pebble adjacent 
to the root and a vertex with at least two pebbles is a {\sf C3} configuration.
Likewise, {\sf S2}, {\sf S3}, and {\sf S4} follow from the nonexistence of {\sf C6},
{\sf C4B}, and {\sf C4A} respectively.
Next, {\sf S5} simply partitions the vertices according to their number of pebbles, 
then uses the definition of $N_0'$.
Finally, since $C$ has no {\sf C5}, no two vertices of $N_3$ are adjacent.
However, because $G$ has diameter two, every such $x$ and $y$ have a common neighbor.
Now the nonexistence of {\sf C7} implies that such common neighbors are distinct, 
which implies {\sf S6}.
\pfc
\medskip

Next we use {\sf S4} and {\sf S5} to count $|C|$ in two ways:
$$3n_3 + 2n_2 + n_1\
=\ n_3 + n_2 + n_1 + (n_{3,r} + n_{3,3} + n_{2,r} + n_0') + 7k - 6.$$
Then {\sf S3} and {\sf S6} imply
\begin{eqnarray*}
0&= & -2n_3 - n_2 + n_{3,r} + n_{3,3} + n_{2,r} + n_0' + 7k - 6\\
&\ge& -n_3 + \binom{n_3}{2} + n_0' + 7k - 6. \\
\end{eqnarray*}
Finally, by completing the square and using $n_0' \ge 1$ (since $r \in N_0'$)
and $k \ge 1$, we have
\begin{eqnarray*}
0&< & (n_3 - 3/2)^2 + (4 - 9/4) \\
& =& 2\left[\binom{n_3}{2} - n_3 + 2\right] \\
&\le& 2\left[\binom{n_3}{2} - n_3 + n_0' + 7k - 6\right] \\
& \le& 0\ ,
\end{eqnarray*}
which is a contradiction.
Hence, $C$ must contain a solution of cost at most 7, afterwhich at least
$n+7(k-1)-6$ pebbles remain, from which we obtain $k-1$ more solutions.
\pf

%
%
\section{Remarks}\label{Remarks}

We believe that the upper bound of Theorem \ref{kpd2} can be tightened by
reducing the coefficient of k.
Doing this requires restricting cheap solutions to lesser cost, which
necessitates considering more of them.
For example, there are one cost-4, one cost-5, and four cost-6 solutions 
that were not used in our argument.
Our lower bound has inspired the next conjecture.

\begin{cnj}\label{d2conj}
If $G$ is a graph on $n$ vertices with $\diam(G)\le 2$ then $\p_k(G)\le n+4k-3$.
\end{cnj}

Of course, the Fractional Pebbling Theorem implies that the coefficient of $k$
is $4$ in the limit; in fact, its proof is based on the pigeonhole principle
--- for large enough $k$, {\sf C4A} exists.
Also, Theorem \ref{kpd2} suggests the following problem.

\begin{prb}\label{diamd}
Find upper bounds for the $k$-pebbling numbers of graphs of diameter $d$.
\end{prb}

Along these lines, only the following result is known, proved by Bukh \cite{B}.

\begin{thm}\label{d3}
If the $diam(G) = 3$, then $\p(G) \le (3/2)n + O(1)$.
\end{thm}

In addition, the following question is still open.

\begin{qst}\label{DiamConn}
Is it possible to lower the connectivity requirement in Result \ref{chkt}?
\end{qst}

The construction in \cite{H} shows that $\kappa\ge 2^d/d$
is necessary.

%
%
\bibliographystyle{plain}
%

%
%
\end{document}